\documentclass[11pt]{amsart}
\usepackage[arrow,matrix]{xy}
\usepackage{amsmath,amssymb, bbm, amscd, amsthm,mathrsfs,hyperref}
\theoremstyle{plain}
\newtheorem{thm}{Theorem}[section]
\newtheorem{cor}[thm]{Corollary}
\newtheorem{lem}[thm]{Lemma}
\newtheorem{prop}[thm]{Proposition}
\theoremstyle{definition}
\newtheorem{defn}[thm]{Definition}
\newtheorem{eg}[thm]{Example}

\newtheorem{rem}[thm]{Remark}

\def\al{\alpha}
\def\bt{\beta}
\def\dt{\delta}

\def\tt{\theta}
\def\vph{\varphi}
\def\lmd{\lambda}
\def\vps{\varepsilon}

\def\bgm{\Gamma}
\def\bd{\Delta}

\def\ms{\mathscr}
\def\mc{\mathcal}
\def\mk{\mathfrak}

\def\ot{\otimes}
\def\op{\oplus}

\def\se{\leqslant}
\def\le{\geqslant}
\def\lan{\langle}
\def\ran{\rangle}

\def\Hom{\operatorname {Hom}}

\def\RHom{\operatorname {RHom}}
\def\Ext{\operatorname {Ext}}

\def\End{\operatorname {End}}
\def\Im{\operatorname {Im}}

\def\ad{\operatorname {ad}}
\def\dim{\operatorname {dim}}
\def\id{\operatorname {id}}

\def\injdim{\operatorname {injdim}}

\def\Gr{\operatorname {Gr}}

\def\hdet{\operatorname {hdet}}

\def\t{\text}

\def\it{\textit}

\def\kk{\mathbbm{k}}
\def\ul{\underline}
\def\ra{\rightarrow}
\def\xra{\xrightarrow}
\def\rap{\rightharpoonup}
\def\ZZ{\mathbb{Z}}

\def\gg{\mathfrak{g}}

\begin{document}
\title{\bf  The Calabi-Yau property of Hopf algebras and braided Hopf algebras}

\author{Xiaolan YU}
\address {Xiaolan YU\newline Department of Mathematics, Hangzhou Normal University, Hangzhou, Zhejiang 310036, China}
 \email{xlyu1005@yahoo.com.cn }

\author{Yinhuo ZHANG}
\address {Yinhuo ZHANG\newline Department WNI, University of Hasselt, Universitaire Campus, 3590 Diepeenbeek,Belgium } \email{yinhuo.zhang@uhasselt.be}

\date{}

\begin{abstract}
Let $H$ be a finite dimensional semisimple Hopf algebra and $R$ a braided Hopf algebra in the category of Yetter-Drinfeld modules over $H$. When $R$ is a Calabi-Yau algebra, a necessary and sufficient condition for $R\#H$ to be a Calabi-Yau Hopf algebra is given. Conversely, when $H$ is the group algebra of a finite group and the smash product $R\#H$ is a Calabi-Yau algebra, we give a necessary and sufficient condition for the algebra $R$ to be a Calabi-Yau algebra.

\end{abstract}


\maketitle

\section*{Introduction}

In recent years, Calabi-Yau (CY) algebras have attracted lots of
attention due to their applications in Algebraic Geometry and in Mathematical Physics.
The study of Calabi-Yau Hopf algebras is initiated by K. Brown and J. Zhang in 2008 cf.\cite{bz}, where they studied rigid dualizing complexes of Noetherian Hopf algebras. S. Chemala showed in \cite{c} that quantum enveloping algebras are Calabi-Yau. In \cite{hvz} J. He, F. Van Oystaeyen and Y. Zhang showed that the smash product of a universal enveloping algebra of a finite dimensional Lie algebra is Calabi-Yau if and only if the group is a subgroup of the special linear group and the enveloping algebra itself is Calabi-Yau. Thus they were able to classify the Noetherian cocommutative Calabi-Yau Hopf algebras of dimension less than 4 over an algebraically closed field. The smash product construction of Calabi-Yau Hopf algebras applied in \cite{hvz} provides in fact an effective method to construct new Calabi-Yau (Hopf) algebras based on existing Calabi-Yau (Hopf) algebras. However, the Calabi-Yau property of the smash product $R\#\kk G$ depends strongly on the action of $\kk G$ on $R$.  For example, the pointed Hopf algebra $U(\mc{D},\lmd)$ of finite Cartan type constructed in \cite{as3} with $\Gamma$ an infinite group of finite rank is Calabi-Yau if and only if the associated graded Hopf algebra $R\#\kk\Gamma$ is Calabi-Yau, where $R$ is the Nichols algebra of $U(\mc{D},\lmd)$. But in this case, if $R\#\kk\Gamma$ is Calabi-Yau, then $R$ can not be Calabi-Yau; and vice versa cf.\cite{yuz}. This arises the question: can we find the "right" action of $G$ on $R$ so that   the Calabi-Yau property of an algebra $R$ delivers the Calabi-Yau property of $R\#\kk G$?

The question was answered by Wu and Zhu in \cite{wz}, where they considered the smash product $R\#\kk G$ of a Koszul Calabi-Yau algebra $R$ by a finite group of automorphisms of $R$. They
showed that the smash product $R\# \kk G$ is Calabi-Yau if and only if the
homological determinant (Definition \ref{hdetl}) of the $G$-action on $R$
is trivial. Later, this result was generalized to the case where $R$
is a $p$-Koszul Calabi-Yau algebra and $\kk G$ is replaced by an involutory Calabi-Yau Hopf algebra \cite{liwz}.

Inspired by the work of Wu and Zhu \cite{wz} and the fact that the associated graded Hopf algebra of a pointed Hopf algebra is a smash product of a braided Hopf algebra in a Yetter-Drinfeld module category over the coradical, we consider in this paper the Calabi-Yau property of a smash product Hopf algebra $R\# H$, where $R$ is a braided Hopf algebra in the Yetter-Drinfeld module category over $H$. We use the homological determinant of the Hopf action to describe the homological integral (Definition \ref{int}) of $R\#H$. We then
give a necessary and sufficient condition for $R\#H$ to be a Calabi-Yau
algebra in case  $R$ is Calabi-Yau and $H$ is semisimple (Theorem \ref{cyrtoh}). We then continue to consider the inverse problem. That is, if $R\#H$ is Calabi-Yau, when is $R$ Calabi-Yau? In Section 3, we answer this question  in case
$H=\kk G$ is the group algebra of a finite group.  We then go on to
characterize the Calabi-Yau property of $R$ when $R\#\kk G$ is Calabi-Yau (Theorem \ref{thm cyator}).  Applying our characterization theorem we obtain the Calabi-Yau property of $U(\mc{D},\lmd)$ in case the datum is not generic (). The generic case is completely worked out in \cite{yuz}. We will provide two examples of Calabi-Yau pointed Hopf algebras with a finite abelian group of group-like elements.

The paper is organized as follows. In Section 1, we review the definition of a braided Hopf algebra, the definition of a Calabi-Yau algebra, the concept of a homological integral and the notion of homological determinants.

In Section 2, we study the Calabi-Yau property under a Hopf action. Our main result in this section is Theorem \ref{cyrtoh}, which characterizes the Calabi-Yau property of the smash product Hopf algebra $R\# H$, where $H$ is a semisimple Hopf algebra and $R$ is a braided Hopf algebra in the Yetter-Drinfeld module category over $H$.

In Section 3, we consider the question when the Calabi-Yau property of $R\# H$ implies that $R$ is Calabi-Yau. We answer this question  in case
$H=\kk G$ is the group algebra of a finite group. We first
construct a bimodule resolution of $R$ from a projective resolution
of $\kk$ over the algebra $R\#\kk G$. Based on this, we obtain a
rigid dualizing complex of $R$ in case $R$ is AS-Gorenstein (Theorem \ref{thm ridfgp}). Our main result in this section is Theorem \ref{thm cyator}.

Throughout, we work over an algebraically closed field of characteristic zero. All vector spaces and algebras are assumed to be over $\kk$.

\section{Preliminaries}

Given an algebra $A$, let $A^{op}$ denote the opposite algebra of
$A$ and $A^e$ denote the enveloping algebra $A\ot A^{op}$ of $A$. The unfurnished tensor $\ot$ means
$\ot_\kk$ in this paper. Mod$A$ denotes the category of left $A$-modules.
We use Mod$A^{op}$ to denote the category of right $A$-modules.

For a left $A$-module $M$ and an algebra automorphism $\phi:A\ra A$,
write $_{\phi}M$ for the left $A$-module defined by $a\cdot m=\phi(a)m$ for any $a\in A$ and $m\in M$.  Similarly, for a right $A$-module $N$, we have
$N_{\phi}$. Observe that  $A_\phi\cong {}_{\phi^{-1}} A$ as
$A$-$A$-bimodules. $A_\phi\cong A$ as $A$-$A$-bimodules if and only
if $\phi$ is an inner automorphism.

For a Hopf algebra, we use  Sweedler's notation for its
comultiplication and its coactions of. Let $A$ be a Hopf algebra, and $\xi:A\ra\kk$ an algebra
homomorphism. We write $[\xi]$ to be the winding homomorphism of
$\xi$ defined by
$$[\xi](a)=\sum\xi(a_1)a_2,$$ for any $a\in A$.

\subsection{Braided Hopf algebra}

Let $H$ be a Hopf algebra. We denote by $^H_H\mc{YD}$ the category of Yetter-Drinfeld modules
over $H$ with  morphisms given by $H$-linear and $H$-colinear maps. If $\bgm$ is a finite group, then $^{\kk\bgm}_{\kk\bgm}\mc{YD}$ will be abbreviated to $^{\bgm}_{\bgm}\mc{YD}$.

Assume that $R$ is a braided Hopf algebra in the
category $^H_H\mc{YD}$. For $h\in H$ and $r\in R$, We write $h(r)$ for $h$ acting on $r$. It is an element in $R$. On the other hand, we write $hr$ for $h$ multiplying with $r$. It is an element in $R\#H$.

The tensor product of two Yetter-Drinfeld modules $M$ and $N$ is
again a Yetter-Drinfeld module with the module and comodule
structures given as follows
$$h(m\ot n)= h_1 m\ot  h_2 n \t{ and } \dt(m\ot  n) =
m_{(-1)}n_{(-1)}\ot  m_{(0)}\ot n_{(0)}, $$for any $h\in H$, $m\in
M$ and $n\in N$.  This turns the category of Yetter-Drinfeld modules
$^H_H\mc{YD}$ into a braided tensor category. For more detail about braided tensor categories, one refers
to \cite{kas}.

For any two Yetter-Drinfeld modules $M$ and $N$, the braiding
$c_{M,N} : M\ot N\ra  N \ot M $ is given by $$ c_{M,N}(m\ot n) =
m_{(-1)}\cdot n \ot m_{(0)},$$ for any $ m \in M$ and $n\in N$.

\begin{defn}Let $H$ be a Hopf algebra.
\begin{enumerate}
\item[(i)] An \textit{algebra} in $^H_H\mc{YD}$ is a $\kk$-algebra $(R,m,u)$ such that
$R$ is a Yetter-Drinfeld $H$-module, and both the multiplication $m:R\otimes R\ra R$ and the unit $u:\kk\ra R$ are morphisms in $^H_H\mc{YD}$.
\item[(ii)] A  \textit{coalgebra}  in $^H_H\mc{YD}$ is
a $\kk$-coalgebra $(C,\Delta,\varepsilon)$ such that $C$ is a Yetter-Drinfeld $H$-module, and both the comultiplication $\bd:R\ra R\otimes R $ and the counit $\vps:R\ra \kk$  are morphisms in
$^H_H\mc{YD}$.
\end{enumerate}
\end{defn}

Let $R$ and $S$ be two algebras in $^H_H\mc{YD}$. Then $R\ot S$ is a
Yetter-Drinfeld module in $^H_H\mc{YD}$, and  becomes an algebra in
the category $^H_H\mc{YD}$ with the  multiplication $m_{R\ul{\ot}S}$
defined by $$m_{R\ul{\ot}S}:=(m_R\ot m_S)(\id\ot c \ot \id).$$
Denote this algebra by $R\ul{\ot} S$.

\begin{defn}\index{braided bialgebra}\index{braided Hopf algebra!graded}\index{braided Hopf algebra}
Let $H$ be a Hopf algebra. A\textit{ braided bialgebra} in
$^H_H\mc{YD}$ is a collection $(R,m, u,\bd, \vps)$, where
\begin{enumerate}
\item[(i)] $(R,m, u)$ is an algebra in $^H_H\mc{YD}$.
\item[(ii)] $(R, \bd, \vps)$ is a coalgebra in $^H_H\mc{YD}$.
\item[(iii)] $\bd:R\ra R\ul\ot R$ and $\vps: R\ra \kk$ are  morphisms of algebras in $^H_H\mc{YD}$.
\end{enumerate}
If, in addition, the identity is convolution invertible in $\End
(R)$, then $R$ is called a \it{braided Hopf algebra} in
$^H_H\mc{YD}$. The inverse of the identity is called the
\textit{antipode} of $R$.

\end{defn}

In order to distinguish comultiplications of braided Hopf algebras
from those of usual Hopf algebras, we use Sweedler's notation with
upper indices for braided Hopf algebras
\begin{equation*}\label{sweed}\bd(r)=r^1\ot r^2.\end{equation*}

Let $H$ be a Hopf algebra and $R$ a braided Hopf algebra in the
category  $^H_H\mc{YD}$. Then $R\#H$ is a usual Hopf algebra with
the following structure \cite{ra}:

The multiplication is given by \begin{equation*}(r\#g)(s\#h):=rg_1(s)\#g_2h\end{equation*}
with unit $u_R\ot u_H$.
The comultiplication is given by
\begin{equation}\label{equa braidedcomulti}
\bd(r\#h):=r^1\#(r^2)_{(-1)}h_1\ot
(r^2)_{(0)}\#h_2
\end{equation}
with counit   $\vps_R\ot \vps_H$.
 The antipode is as follows:
\begin{equation}\label{equa braidedanti}\mc{S}_{R\#H}(r\#h)=(1\#\mc{S}_H(r_{(-1)}h))(\mc{S}_R(r_{(0)})\#1).
\end{equation}

The algebra $R\#H$ is called the \textit{Radford biproduct} or
\textit{bosonization} of $R$ by $H$. The algebra $R$ is a subalgebra of $R\#H$
and $H$ is a Hopf subalgebra of $R\#H$. \index{Radford biproduct}
\index{bosonization}

Conversely, let $A$ and $H$ be two Hopf algebras and $\pi:A\ra H$,
$\iota:H\ra A$ Hopf algebra homomorphisms such that
$\pi\iota=\id_H$. In this case the algebra of right coinvariants
with respect to $\pi$
$$R=A^{co\pi}:=\{a\in A\mid (\id\ot \pi)\bd(a)=a\ot 1\},$$
is a braided Hopf algebra in $^H_H\mc{YD}$, with the following
structure \cite{ra}:
\begin{itemize}
\item[(i)] The action of $H$ on $R$ is the restriction of the adjoint action (composed with $\iota$).
\item[(ii)] The coaction is $(\pi\ot \id)\bd$.
\item[(iii)] $R$ is a subalgebra of $A$.
\item[(iv)] The comultiplication is given by  $$\bd_R(r)=r_1\iota\mc{S}_H\pi (r_2)\ot
r_3.$$
\item[(v)] The antipode is given by $$\mc{S}_R(r)=\pi(r_1)\mc{S}_A(r_2).$$
\end{itemize}

Define a linear map $\rho: A\ra R$ by
$$\rho(a)=a_1\iota\mc{S}_H\pi(a_2),$$ for all $r\in R$.

\begin{thm}\cite{ra}
The morphisms $\Psi:A\ra R\#H$ and $\Phi: R\#H\ra A$ defined by
$$\Psi(a)=\rho(a_1)\#\pi(a_2)\t{ and } \Phi(r\#h)=r\iota(h)$$ are mutually inverse isomorphisms of Hopf algebras.
\end{thm}

\subsection{Calabi-Yau algebras} We follow Ginzburg's definition of a Calabi-Yau algebra \cite{g2}.

\begin{defn}An algebra $A$ is called a \textit{Calabi-Yau algebra of dimension
$d$} if
\begin{enumerate}
\item[(i)] $A$ is \it{homologically smooth}, that is, $A$ has
a bounded resolution of finitely generated projective
$A$-$A$-bimodules;
\item[(ii)] There are $A$-$A$-bimodule
isomorphisms$$\Ext_{A^e}^i(A,A^e)=\begin{cases}0& i\neq d
\\A&i=d.\end{cases}$$
\end{enumerate}
\end{defn}
In the sequel, Calabi-Yau will be abbreviated to CY for short.

CY algebras form a class of algebras possessing a rigid dualizing complex. The
non-commutative version of a dualizing complex was first introduced
by Yekutieli.

A Noetherian algebra in this paper means a \textit{left and right}
Noetherian algebra.

\begin{defn}\cite{y} (cf. \cite[Defn. 6.1]{vdb})\label{defn dc} \index{dualizing complex}
Assume that $A$ is a  (graded) Noetherian algebra. Then an object
$\ms{R}$ of $D^b(A^e)$ ($D^b(\t{GrMod}(A^e))$) is called a
\textit{dualizing complex} (in the graded sense) if it satisfies the
following conditions:
\begin{enumerate}
\item[(i)] $\ms{R}$ is of  finite injective dimension over $A$ and $A^{op}$.
\item[(ii)] The cohomology of $\ms{R}$ is given by bimodules which are
finitely generated on both sides.
\item[(iii)] The natural morphisms $A\ra  {\RHom}_A(\ms{R},\ms{R})$ and $A\ra
\RHom_{A^{op}}(\ms{R},\ms{R})$ are isomorphisms in $D(A^e)$
($D(\t{GrMod}(A^e))$).
\end{enumerate}
\end{defn}

Roughly speaking, a dualizing complex is a complex $\ms{R}\in
D^b(A^e)$ such that the functor
\begin{equation}\label{dualcom}\RHom_A(-,\ms{R}):D^b_{fg}(A)\ra D^b_{fg}(A^{op})\end{equation} is a
duality, with adjoint $\RHom_{A^{op}}(-,\ms{R})$ (cf. \cite[Prop.
3.4 and Prop. 3.5]{y}). Here $D^b_{fg}(A)$ is the  full triangulated
subcategory of $D(A)$ consisting of bounded complexes with finitely
generated cohomology modules.

In the above definition,  the algebra $A$ is a Noetherian algebra.
In this case, a dualizing complex in the graded sense is also a
dualizing complex in the usual sense.

Dualizing complexes are not unique  up to isomorphism. To overcome
this weakness,  Van den Bergh introduced the concept of a rigid
dualizing complex cf. \cite[Defn. 8.1]{vdb}.

\begin{defn}\label{defn rdc}\index{dualizing complex!rigid} Let $A$ be a (graded) Noetherian algebra. A dualizing
 complex $\ms{R}$ over $A$ is called \textit{rigid} (in the graded sense) if $$\RHom_{A^e}(A,{_A\ms{R}\ot \ms{R}_A})\cong
 \ms{R}$$ in $D(A^e)$ ($D(\t{GrMod}(A^e))$).
\end{defn}

Note that if $A^e$ is Noetherian, then the graded version of
this definition implies the ungraded version. The following lemma can be found in  \cite[Prop. \;4.3]{bz}  and  \cite[Prop. \;8.4]{vdb}.
Note that if a Noetherian algebra has finite left and
right injective dimension, then they are equal cf. \cite[Lemma
A]{za}. We call this common value the injective dimension of $A$.

\begin{lem}\label{vdb} Let $A$ be a Noetherian algebra.   Then the following two conditions are
equivalent:
\begin{enumerate}
\item[(i)] $A$ has a rigid dualizing complex $\ms{R}=A_{\psi}[s]$, where $\psi$ is an
algebra automorphism and $s\in \ZZ$.
\item[(ii)] $A$ has finite
injective dimension $d$ and there is an algebra automorphism $\phi$
such that
$$\Ext_{A^e}^i(A,A^e)\cong\begin{cases}0,& i\neq d;
\\A_\phi&i=d\end{cases}$$ as $A$-$A$-bimodules. \end{enumerate}
If the two conditions are equivalent, then
 $\phi=\psi^{-1}$ and $s=d$.
\end{lem}

 The
following corollary follows directly from Lemma \ref{vdb} and the
definition of a CY algebra.

\begin{cor}\label{cor cyrid}
Let $A$ be a Noetherian algebra which is  homologically smooth. Then
$A$ is a CY algebra of dimension $d$ if and only if $A$ has a rigid
dualizing complex $A[d]$.
\end{cor}

\subsection{Homological integral}

In \cite{hvz}, the CY property of Hopf algebras was discussed by
using the homological integrals of Artin-Schelter Gorenstein (AS-Gorenstein for short)
algebras \cite[Thm. 2.3]{hvz}. The concept of a homological integral for an
AS-Gorenstein Hopf algebra was introduced by Lu,
Wu and Zhang in \cite{lwz} to study infinite dimensional Noetherian Hopf algebras. It generalizes
the concept of an integral of a finite dimensional Hopf algebra. In \cite{bz},
homological integrals were defined for general AS-Gorenstein
algebras.

\begin{defn}(cf. \cite[defn. 1.2]{bz}). \label{defn as}\index{AS-regular}
\begin{enumerate}
\item[(i)] Let $A$ be a  left Noetherian augmented algebra with
a fixed augmentation map $\vps:A\ra \kk$. $A$ is said to be \it{left
AS-Gorenstein},\index{AS-Gorenstein!left} if
\begin{enumerate}
\item[(a)] $\injdim  {_AA}=d<\infty$,
\item[(b)] $\dim\Ext_A^i({_A\kk},{_AA})=\begin{cases}0,&i\neq d;\\ 1,& i=d,\end{cases}$
\end{enumerate}
where $\injdim$ stands for  injective dimension.

A \textit{Right AS-Gorenstein algebras}\index{AS-Gorenstein!right} can be defined similarly.

\item[(ii)] An algebra $A$ is
said to be \textit{AS-Gorenstein}\index{AS-Gorenstein} if it is both left and right
AS-Gorenstein (relative to the same augmentation map $\vps$).
\item[(iii)] An
AS-Gorenstein algebra $A$ is said to be \it{regular}  if in
addition, the global dimension of $A$
 is finite.
\end{enumerate}
\end{defn}


\begin{defn}\label{int}
 Let $A$ be a left AS-Gorenstein algebra
with $\injdim{_AA}=d$. Then $\Ext_A^d({_A\kk},{_AA})$ is a 1-dimensional
right $A$-module. Any nonzero element in $\Ext_A^d({_A\kk},{_AA})$
is called a \it{left homological integral} of $A$. We write
$\int_A^l$ for $\Ext_A^d({_A\kk},{_AA})$. Similarly, if $A$ is right
AS-Gorenstein with $\injdim{A_A}=d$, any nonzero element  in $\Ext_A^d({\kk_A},{A_A})$ is
called a \it{right homological integral} of $A$. Write
$\int_A^r$ for $\Ext_A^d({\kk_A},{A_A})$.

$\int_A^l$ and $\int_A^r$
are  called \textit{left and right homological integral modules} of $A$
respectively.
\end{defn}

The left integral module $\int^l_A$ is a 1-dimensional right
$A$-module. Thus $\int^l_A\cong \kk_\xi$ for some algebra
homomorphism $\xi:A\ra\kk$.

\begin{prop}\label{prop cy-as}
Let $A$ be a Noetherian augmented algebra such that $A$ is CY of
dimension $d$. Then $A$ is AS-regular of global dimension $d$. In
addition, $\int_A^l\cong \kk$ as right $A$-modules.
\end{prop}
\proof If $A$ is an augmented algebra, then $_A\kk$ is a finite
dimensional $A$-module. By \cite[Remark 2.8]{br}, $A$ has global
dimension $d$.

It follows from  \cite[Prop. 2.2]{br} that $A$ admits a projective
bimodule resolution
$$0\ra P_d\ra \cdots\ra P_1\ra P_0\ra A\ra 0,$$ where each $P_i$ is
finitely generated as an $A$-$A$-bimodule. Tensoring with functor
$\ot_A \kk$, we obtain a projective resolution of
 $_A\kk$:
$$0\ra
P_d\ot_A\kk\ra \cdots\ra P_1\ot_A\kk\ra P_0\ot_A\kk\ra {}_A\kk\ra
0.$$ Since each $P_i$ is finitely generated, the isomorphism
$$\kk\ot _A\Hom_{A^e}(P_i,A^e)\cong \Hom_{A }(P_i\ot_A\kk,A)$$ holds
in Mod$A^{op}$. Therefore, the complex $\Hom_A(P_\bullet\ot_A\kk,
A)$ is isomorphic to the complex
$\kk\ot_A\Hom_{A^e}(P_\bullet,A^e)$. The fact that algebra $A$ is CY of
dimension $d$ implies that the following $A$-$A$-bimodule complex is exact:
$$0\ra \Hom_{A^e}(P_0,A^e)\ra \cdots\ra \Hom_{A^e}(P_{d-1},A^e)\ra  \Hom_{A^e}(P_d,A^e)\ra A\ra 0.$$ Thus the complex $\kk\ot_A\Hom_{A^e}(P_\bullet,A^e)$ is exact except at $\kk\ot_A\Hom_{A^e}(P_d,A^e)$, whose homology is $\kk$.
It follows that the isomorphisms
$$\Ext_A^i({}_A\kk,{}_AA)\cong\begin{cases}0,&i\neq d;\\\kk,&i=
d\end{cases}$$ hold in Mod$A^{op}$. Similarly, we have isomorphisms
$$\Ext_A^i( \kk_A, A_A)\cong\begin{cases}0,&i\neq d;\\\kk,&i=
d\end{cases}$$ in Mod$A$. We conclude that $A$ is AS-regular and
$\int^l_A\cong \kk$. \qed

\begin{rem}\label{twisted}
From the proof of Proposition \ref{prop cy-as} we can see that if
$A$ is a Noetherian augmented algebra such that
\begin{enumerate}\item[(i)] $A$ is  homologically
smooth, and
\item[(ii)] there is an integer $d$ and an algebra automorphism $\psi$, such that
$$\Ext_{A^e}^i(A,A^e)\cong\begin{cases}0,& i\neq d;
\\A_\psi,&i=d\end{cases}$$ as $A$-$A$-bimodules,
\end{enumerate}
then $A$ is  AS-regular of global dimension $d$. In this case,
$\int_A^l\cong \kk_\xi$. The algebra homomorphism $\xi$ is defined
by $\xi(a)=\vps(\psi(a))$ for all $a\in A$, where $\vps$ is the
augmentation map of $A$.
\end{rem}

\subsection{Homological determinants}The homological determinant  was defined by J{\o}rgensen and Zhang
\cite{joz} for graded automorphisms of an
AS-Gorenstein algebra and by Kirkman, Kuzmanovich and
Zhang \cite{kkz} for Hopf actions on an AS-Gorenstein algebra. The homological determinant was used to study the
AS-Gorenstein property of invariant subrings.

\begin{defn}\label{hdetl}(cf. \cite{liwz}, \cite{kkz})
Let $H$ be a Hopf algebra, and  $R$ an $H$-module AS-Gorenstein algebra of injective
dimension $d$.  There is a natural  $H$-action on $\Ext_R^d(\kk,R)$ induced by the
$H$-action on $R$. Let $\bf{e}$ be a non-zero element in
$\Ext_R^d(\kk,R)$. Then there exists an algebra homomorphism $\eta:
H\ra\kk$ satisfying $h \cdot {\bf e} = \eta(h)\bf{e}$ for all $h\in
H$.
\begin{enumerate}
\item[(i)] The
composite map $\eta\mc{S}_H:H\ra\kk$  is called the \textit{
homological determinant} of the $H$-action on $R$, and it is denoted
by $\hdet $ (or more precisely $\hdet_R$).
\item[(ii)] The homological determinant $\hdet_R$ is said to be \textit{trivial} if $\hdet_R=\vps_H$, where $\vps_H$ is the counit of the Hopf algebra $H$.
\end{enumerate}
\end{defn}\index{homological determinant}

\section{Calabi-Yau property under Hopf actions}\label{fgp. H}

Let $H$ be a Hopf algebra and $R$ a braided Hopf algebra in the
category $^H_H\mc{YD}$. In this section, we study the CY property of the smash product $R\#H$, when $R$ is a CY algebra and $H$ is a  semisimple Hopf algebra.

 For a left $R\#H$-module $M$, the vector
space $M\ot H$ is a left $R\#H$-module defined by
$$(r\#h)\cdot (m\ot g):=(r\#h_1)m\ot h_2g,$$ for  all $r\#h\in R\#H$ and $m\ot g\in M\ot H$. Denote this $R\#H$-module by $M\#H$.

Let $M$ and $N$ be two $R\#H$-modules. Then there is a natural left
$H$-module structure on $\Hom_R(M,N)$ given by the adjoint action
$$(h\rap f)(m):=h_2f(\mc{S}^{-1}_H(h_1)m),$$ for all $h\in H$, $f\in \Hom_R(M,N)$ and $m\in M$.

\begin{lem}\label{lem H1}
Let $M$ be a left $R\#H$-module. Then $\Hom_R(M,R)\ot H$ is an
$H$-$R\#H$-bimodule, where the left $H$-module structure is defined
by
$$h\cdot (f\ot g):=h_1\rap f\ot h_2g$$ and the right $R\#H$-module structure is given by
the diagonal action:
$$(f\ot g)\cdot (r\#h):=f  g_1(r)\ot g_2h,$$ for all $f\in \Hom_R(M,R)$, $g,h\in H$ and $r\in R$.
\end{lem}
\proof First we show that for all $h\in H$, $f\in \Hom_R(M,R)$ and
$r\in R$\begin{equation}\label{41}(h_1\rap f){h_2}(r)=h
\rightharpoonup(fr).\end{equation}
For $m\in M$, we have $$\begin{array}{ccl}[(h_1\rap f){h_2}(r)](m)&=&(h_1\rap f)(m)h_2(r)\\
&=&h_2(f(\mc{S}^{-1}_H(h_1)m))h_3(r)\\
&=&h_2(f(\mc{S}^{-1}_H(h_1)m)r)\\
&=&h_2((fr)(\mc{S}^{-1}_H(h_1)m))\\
&=&[h\rap(fr)](m).
\end{array}$$

Now we check that for  all $f\ot g\in \Hom_R(M,R)\ot H$, $h\in H$
and $r\#k\in R\#H$, $(h\cdot(f\ot g))\cdot (r\# k)=h\cdot((f\ot
g)\cdot (r\# k))$. We have
$$\begin{array}{ccl}(h\cdot(f\ot g))\cdot (r\# k)&=&(h_1\rap f\ot h_2g)\cdot (r\# k)\\
&=&(h_1\rap f) (h_2g_1)(r)\ot h_3g_2k.
\end{array}$$ and
$$\begin{array}{ccl}h\cdot((f\ot g)\cdot (r\# k))&=&h\rap(fg_1(r)\ot g_2k)\\
&=&h_1\rap (fg_1(r))\ot h_2g_2k\\
&\overset{\tiny{(\ref{41})}}=&(h_1\rap f)(h_2g_1)(r)\ot h_3g_2k.
\end{array}$$\qed

There is a natural right $R\#H$-module structure on $\Hom_{R\#H}(M\#
H,R\#H)$. It is  also a left $H$-module defined by
\begin{equation}\label{H2}(h\cdot f)(m\ot g):=f(m\ot gh),\end{equation} for all $h\in H$, $f\in
\Hom_{R\#H}(M\# H,R\#H)$ and $m\ot g\in M\ot H$. Then $\Hom_{R\#H}(M\ot H,R\#H)$ is an $H$-$R\#H$-bimodule.


\begin{prop}\label{prop iso}
Let  $P$ be an $R\#H$-module such that it is finitely generated
projective as an $R$-module. Then
$$\Hom_R(P,R)\ot H\cong \Hom_{R\#H}(P\# H, R\#H)$$ as
$H$-$R\#H$-bimodules.
\end{prop}
\proof Let $$\psi:\Hom_R(P,R)\ot H\ra \Hom_{R\#H}(P\# H, R\#H)$$ be
the homomorphism  defined by
$$\begin{array}{ccl}[\psi(f\ot h)](p\ot g)&=&(g_1\rightharpoonup f)(p)\#g_2h\\
&=&g_2(f(\mc{S}_H^{-1}(g_1)p))\#g_3h,
\end{array}$$ for all $f\ot h\in \Hom_R(P,R)\ot H$ and $p\ot g\in P\# H$.

We claim that the image of $\psi$ is contained in $\Hom_{R\#H}(P\#
H, R\#H)$. For any $f\ot h\in \Hom_R(P,R)\ot H$, $r\#k\in R\#H$ and
$p\ot g\in P\# H$, on one hand, we have
$$\begin{array}{ccl}[\psi(f\ot h)]((r\#k)(p\ot g))&=&[\psi(f\ot h)]((r\#k_1)p\ot k_2g))\\
&=&(k_3g_2)( f(\mc{S}_H^{-1}(k_2g_1)((r\#k_1)p)))\#k_4g_3h\\
&=&(k_2g_3)(
f(((\mc{S}_H^{-1}(k_1g_2))(r))\mc{S}_H^{-1}(g_1)p))\#k_3g_4h.
\end{array}$$
On the other hand,  $$\begin{array}{ccl}(r\#k)[\psi(f\ot h)](p\ot g)&=&(r\#k)(g_2 (f(\mc{S}_H^{-1}(g_1)p))\#g_3h)\\
&=&r(k_1g_2) (f(\mc{S}_H^{-1}(g_1)p))\#k_2g_3h\\
&=&(k_2g_3) (\mc{S}_H^{-1}(k_1g_2)(r) f(\mc{S}_H^{-1}(g_1)p))\#k_3g_4h\\
&=&(k_2g_3)(
f(((\mc{S}_H^{-1}(k_1g_2))(r))\mc{S}_H^{-1}(g_1)p))\#k_3g_4h.
\end{array}$$
 Now we show that $\psi$ is an $H$-$R\#H$-bimodule
homomorphism. We have
$$\begin{array}{ccl}[\psi((f\ot h)(r\#k))](p\ot g)&=&[\psi(fh_1(r)\ot h_2k))](p\ot g)\\
&=&g_2([fh_1(r)](\mc{S}_H^{-1}(g_1)p))\ot g_3h_2k\\
&=&g_2(f(\mc{S}_H^{-1}(g_1)p))(g_3h_1)(r)\ot g_4h_2k\\
&=&(g_2(f(\mc{S}_H^{-1}(g_1)p))\ot g_3h)(r\# k)\\
&=&[\psi(f\ot h)(r\#k)](p\ot g)
\end{array}$$
 and
$$\begin{array}{ccl}[\psi(k(f\ot h))](p\ot g)&=&[\psi(k_1\rap f\ot k_2h)](p\ot g)\\
&=&g_2((k_1\rap f)(\mc{S}_H^{-1}(g_1)p))\#g_3k_2h\\
&=&(g_2k_2)(f(\mc{S}_H^{-1}(k_1)\mc{S}_H^{-1}(g_1)p))\#g_3k_3h\\
&=&((g_1k_1)\rap f)(p)\ot g_2k_2h\\
&=&[\psi(f\ot h)](p\ot gk)\\
&=&[k\cdot \psi(f\ot h)](p\ot g).
\end{array}$$
So $\Hom_R(P,R)\ot H\cong \Hom_{R\#H}(P\# H, R\#H)$ as
$H$-$R\#H$-bimodules when $P$ is finitely generated projective as an
$R$-module. \qed

\begin{prop}\label{ext} Let $H$ be a finite dimensional Hopf algebra and $R$ a Noetherian braided Hopf algebra in the category $^H_H\mc{YD}$. Then
$$\Ext^i_{R\#H}(H,R\#H)\cong \Ext^i_R(\kk,R)\ot H$$ as $H$-$R\#H$-bimodules for all $i\le 0$.
\end{prop}
\proof Since $R$ is Noetherian and $H$ is finite dimensional, $R\#H$
is also Noetherian. Then $_{R\#H}\kk$ admits a projective resolution
 $$\cdots\ra P_n\ra \cdots\ra P_1\ra P_0\ra \kk\ra 0$$
such that each  $P_n$ is a finitely generated $R\#H$-module. Because
$H$ is finite dimensional, each $P_n$ is also finitely generated as
an $R$-module. Tensoring with $H$, we obtain a projective resolution
of $H$ over $R\#H$
$$\cdots\ra P_n\# H\ra \cdots\ra P_1\# H\ra P_0\# H\ra H\ra 0.$$  Applying the functor $\Hom_{R\#H}(-,R\#H)$ to this complex, we obtain the following complex
\begin{equation}\label{comp1}0\ra \Hom_{R\#H}(P_0\# H,R\#H)\ra \Hom_{R\#H}(P_1\# H,R\#H)\ra\cdots \end{equation}
$$\hspace{65mm}\ra \Hom_{R\#H}(P_n\# H,R\#H)\ra \cdots.$$ This is a complex of $H$-$R\#H$-bimodules, where the left $H$-module structure is defined as in (\ref{H2}). By Proposition \ref{prop iso}, one can check that it is isomorphic to the following complex of $H$-$R\#H$-bimodules,
\begin{equation}\label{comp2}0\ra \Hom_R(P_0,R)\ot H\ra \Hom_R(P_1,R)\ot H\cdots \hspace{20mm}\end{equation}
$$\hspace{60mm}\ra \Hom_R(P_n,R)\ot H\ra \cdots.$$
After taking cohomologies of complex (\ref{comp1}) and complex
(\ref{comp2}),  we arrive at isomorphisms of $H$-$R\#H$-bimodules
$$\Ext^i_{R\#H}(H,R\#H)\cong \Ext^i_R(\kk,R)\ot H$$ for all $i\le 0$.\qed

The algebra $R$ can be viewed as an augmented  right $H$-module
algebra through the right $H$-action: $r\cdot
h:=\mc{S}_H^{-1}(h)\cdot r$, for all $r\in R$ and $h\in H$.   The
algebra $H\#R$ can be defined in a similar way. The multiplication
is given by
$$(h\#s)(k\#r):=hk_2\#(s\cdot k_1)r=hk_2\#(\mc{S}_H^{-1}(k_1)(s))r,$$
for all
 $h\#s$ and $k\#r\in H\#R$. The homomorphism $\vph:R\#H\ra H\#R$
defined by
$$\vph(r\#k)=k_2\#\mc{S}_H^{-1}(k_1)(r)$$ is an algebra isomorphism
with its inverse $\psi:H\#R\ra R\#H$ defined by
$$\psi(k\#r)=k_1(r)\#k_2.$$ In addition, $\vph$ is compatible with the
augmentation maps of $R\#H$ and $H\#R$ respectively. Now right
$R\#H$-modules can be treated as $H\#R$-modules. Let $M$ and  $N$ be
two $H\#R$-modules, then $\Hom_R(M,N)$ is a right $H$-module defined
by
$$(f\leftharpoonup h)(m):=f(m\mc{S}_H(h_1))h_2,$$ for all $h\in H$, $f\in \Hom_R(M,N)$ and $m\in M$.

Similar to the left case, we have the following proposition.

\begin{prop}\label{extright}
Let $H$ be a finite dimensional Hopf algebra and $R$ a Noetherian
braided Hopf algebra in the category $^H_H\mc{YD}$. Then
$$\Ext^i_{R\#H}(H_{R\#H},{R\#H}_{R\#H})\cong H\ot
\Ext^i_R(\kk_{R},R_R)$$ as $R\#H$-$H$-bimodules for all $i\le 0$.
\end{prop}

\begin{lem}\label{lgl}
Let $H$ be a Hopf algebra and $R$ an $H$-module algebra.  If the
left global dimensions of $R$ and $H$ are $d_R$ and $d_H$
respectively, then the  left global dimension of $A=R\#H$ is not
greater than $d_R+d_H$.
\end{lem}
\proof Let $M$ and $N$ be two $A$-modules. We have
$$\Hom_A(M,N)\cong \Hom_H(\kk,\Hom_R(M,N)),$$ that is, the functor
$\Hom_A(M,-)$ factors through as follows
$$\xymatrix{
  \t{Mod} A\ar[rr]^{\Hom_R(M,-)} \ar[dr]_{\Hom_A(M,-)}
                &  &    \t{Mod}H \ar[dl]^{\Hom_H(\kk,-)}    \\
                & \t{Mod}\kk              }.$$
To apply the Grothendieck spectral sequence (see e.g. \cite[Sec.
5.8]{we}), we need to show that if $N$ is an injective $A$-module,
then $\Ext_H^q(\kk,\Hom_R(M,N))=0$ for all $q\le 1$.

Let $$\cdots\ra P_i\ra P_{i-1}\ra \cdots \ra P_1\ra P_0\ra \kk\ra
0$$ be a projective resolution of $\kk$ over $H$.
$\Ext_H^*(\kk,\Hom_R(M,N))$ are the cohomologies of the complex
$\Hom_H(P_\bullet, \Hom_R(M,N))$. The following
isomorphisms hold:
$$\begin{array}{ccl}\Hom_H(P_\bullet, \Hom_R(M,N))&\cong&\Hom_H(\kk,\Hom_\kk(P_\bullet,
\Hom_R(M,N)))\\
&\cong &\Hom_H(\kk,\Hom_R(P_\bullet\ot M,
N))\\
&\cong &\Hom_{R\#H}(P_\bullet\ot M, N).
\end{array}$$
 The complex $P_\bullet$ is exact except at $P_0$. Since the functors $\Hom_{R\#H}(-,N)$ and $-\ot M$ are exact, the complex $\Hom_H(P_\bullet, \Hom_R(M,N))$ is also exact except at $\Hom_H(P_0, \Hom_R(M,N))$.  It
follows  that $$\Ext_H^q(\kk,\Hom_R(M,N))=0$$ for all $q\le 1$.

Now we have  $$\Ext^q_H(\kk,\Ext_R^p(M,N))\Rightarrow
\Ext^{p+q}_{R\#H}(M,N).$$ Because the left  global dimensions of $R$
and $H$ are $d_R$ and $d_H$, $\Ext^i_{R\#H}(M,N)=0$ for all $i\le
d_R+d_H$. Therefore, the left global dimension of $R\#H$ is not
greater than $d_R+d_H$.\qed

Let $H$ be an involutory CY Hopf algebra and  $R$   a $p$-Koszul CY
algebra and a left  $H$-module algebra.   As mentioned in
the introduction, Liu, Wu and Zhu used
the homological determinant of the $H$-action to characterize  the
CY property of $R\#H$ in \cite{liwz}. They defined an $H$-module structure on the
Koszul bimodule complex of $R$ and computed the $H$-module
structures on the Hochschild cohomologies. Then they proved that
$R\#H$ is CY if and only if the homological determinant is trivial. If $H$ is not involutory or $R$ is not a $p$-Koszul algebra, is it still true that $R\#H$ is a CY algebra when the homological determinant is trivial?

We answer this question in the case that $R$ is a braided Hopf algebra in the
category $^H_H\mc{YD}$, where $H$ is a semisimple Hopf
algebra.  We use the homological determinant to discuss the
homological integral and the rigid dualizing complex of the algebra
$A=R\#H$. We then give a necessary and sufficient condition for $A$
to be a CY algebra. The result we will obtain is slightly different from what
was obtained by Liu, Wu and Zhu. We first need the following lemma.

\begin{lem}\label{s^2}
Let $H$ be a Hopf algebra, and $R$  a  braided Hopf algebra in the
category ${^H_H\mathcal{YD}}$. Then
$$\mc{S}_{R\#H}^2(r)={\mc{S}_H(r_{(-1)})}({ \mc{S}_R^2(r_{(0)})}),$$
for any $r\in R$.
\end{lem}
\proof  Set $A=R\#H$. By equation (\ref{equa braidedanti}), for any
$r\in R$, we have
$$\mc{S}_A(r)=(1\#\mc{S}_H(r_{(-1)}))(\mc{S}_R(r_{ (0) })\#1).$$
Therefore,
$$\begin{array}{ccl}\mc{S}_A^2(r)&=&\mc{S}_A((1\#\mc{S}_H(r_{(-1)}))(\mc{S}_R(r_{(0)})\#1))\\
&=&\mc{S}_A(\mc{S}_R(r_{(0)})\#1)\mc{S}_A(1\#\mc{S}_H(r_{(-1)}))\\
&=&(1\#\mc{S}_H(\mc{S}_R(r_{(0)})_{(-1)}))(\mc{S}_R(\mc{S}_R(r_{(0)})_{(0)})\#1)(1\#\mc{S}_H^2(r_{(-1)}))\\
&=&(1\#\mc{S}_H(r_{(0)(-1)}))(\mc{S}_R^2(r_{(0)(0)})\#1)(1\#\mc{S}_H^2(r_{(-1)}))\\
&=&(1\#\mc{S}_H(r_{(-1)2}))(\mc{S}_R^2(r_{(0)})\#1)(1\#\mc{S}_H^2(r_{(-1)1}))\\
&=&{\mc{S}_H(r_{(-1)3})}( \mc{S}_R^2(r_{(0)}))\#\mc{S}_H(r_{(-1)2})\mc{S}_H^2(r_{(-1)1})\\
&=&{{\mc{S}_H(r_{(-1)2})} ( \mc{S}_R^2(r_{(0)}))}\#\mc{S}_H(\varepsilon(r_{(-1)1}))\\
&=&{{\mc{S}_H(r_{(-1)})}( \mc{S}_R^2(r_{(0)}))}.\end{array}$$\qed

\begin{prop}\label{thm int}
Let $H$ be a finite dimensional semisimple Hopf algebra and  $R$ a
braided Hopf algebra in the category $^H_H\mc{YD}$. If  $R$ is an
 AS-regular algebra of  global dimension $d_R$, then $A=R\#H$ is  also AS-regular of global dimension
$d_R$.

In this case, $\int^l_A=\kk_\xi$, where $\xi:A\ra \kk$ is defined by
$$\xi(r\#h)=\xi_R(r)\hdet(h),$$ for all $r\#h\in R\#H$,  where the algebra map $\xi_R:R\ra \kk$ defines the left integral module of $R$, i.e., $\int^l_R=\kk_{\xi_R}$.
The rigid dualizing complex of $A$ is isomorphic to $_\psi A[d_R]$,
where $\psi$ is the algebra automorphism $[\xi]\mc{S}_A^2$. To be
precise, $\psi$ is defined by {
$$\psi(r\#h)=\xi_R(r^1)\hdet((r^2)_{(-1)1}h_1)\mc{S}_H((r^2)_{(-1)2})(\mc{S}_R^2((r^2)_{(0)}))\#\mc{S}^2_H(h_2),$$}
for  all $r\#h\in R\#H$.
\end{prop}
\proof  We have the
following isomorphisms
$$\begin{array}{ccl}\RHom_A(\kk,{}_AA)&\cong&\RHom_A(H\ot_{H}\kk,A)\\
&\cong&\RHom_{H}(\kk,\RHom_A(H,A))\\
&\cong&\RHom_{H}(\kk,H)\ot^L_{H}\RHom_A(H,A).
\end{array}$$
Following Proposition \ref{ext}, we have $\int_A^l\cong
\int_H^l\ot_H\int_R^l \ot H$ and
$$\dim\Ext_A^i(\kk,{}_AA)=\begin{cases}0,&i\neq d_R;\\1,&i=
d_R.\end{cases}$$  Let ${\bf e}$ be a non-zero element in $\int_R^l$
and ${\bf h}$ a non-zero element in $\int_H^l$. Since $H$ is  semisimple, $H$ is unimodular. That is, we have $\int_H^l=\kk$. Let $\eta:H\ra\kk$
be an algebra homomorphism such that $h \cdot {\bf e} =
\eta(h)\bf{e}$ for all $h\in H$. Then the following equations hold
$$\begin{array}{ccl}({\bf h}\ot {\bf e}\ot 1)\cdot(r\#h)&=&\xi_R(r){\bf h}\ot {\bf e}\ot h \\
&=&\xi_R(r){\bf h}\ot \vps(h_1){\bf e}\ot h_2 \\
&=&\xi_R(r){\bf h}\ot \eta(\mc{S}_H(h_1))\eta(h_2){\bf e}\ot h_3 \\
&=&\xi_R(r)\eta(\mc{S}_H(h_1)) {\bf h}\ot h_2\cdot({\bf e}\ot 1) \\
&=&\xi_R(r)\eta(\mc{S}_H(h_1)) \varepsilon(h_2){\bf h}\ot {\bf e}\ot 1 \\
&=&\xi_R(r)\hdet(h) {\bf h}\ot {\bf e}\ot 1.
\end{array}$$ This implies that $\int^l_A\cong \kk_\xi$, where $\xi$ is the algebra
homomorphism defined in the proposition.  Similarly, by Proposition
\ref{extright}, we have
$$\dim\Ext_A^i(\kk,A_A)=\begin{cases}0,&i\neq d_R,\\1,&i=
d_R.\end{cases}$$

Because $H$ is finite dimensional and $R$ is Noetherian, the algebra
$A$ is Noetherian as well. Therefore, the left and right global
dimensions of $A$ are equal. Since $H$ is semisimple, the global
dimension of $H$ is 0.  Now it follows from Lemma \ref{lgl} that the
global dimension of $A$ is $d_R$. In conclusion, we have proved that
$A$ is an AS-regular algebra.

By \cite[Prop. \;4.5]{bz}, the rigid dualizing complex of $A$ is
isomorphic to $_{[\xi]\mc{S}_A^2}A[d_R]$. For any $r\#h\in R\#H$, we have
{\small $$\begin{array}{cl}&[\xi]\mc{S}_A^2(r\#h)\\\overset{(a)}=&\mc{S}_A^2[\xi](r\#h)\\\overset{(b)}=&\xi(r^1\#(r^2)_{(-1)}h_1)\mc{S}_A^2((r^2)_{(0)}\#h_2)\\
=&\xi_R(r^1)\hdet((r^2)_{(-1)}h_1)\mc{S}_A^2((r^2)_{(0)})\#\mc{S}^2_H(h_2)\\
\overset{(c)}=&\xi_R(r^1)\hdet((r^2)_{(-1)}h_1)\mc{S}_H((r^2)_{(0)(-1)})(\mc{S}_R^2((r^2)_{(0)(0)}))\#\mc{S}^2_H(h_2)\\
=&\xi_R(r^1)\hdet((r^2)_{(-1)1}h_1)\mc{S}_H((r^2)_{(-1)2})(\mc{S}_R^2((r^2)_{(0)}))\#\mc{S}^2_H(h_2).
\end{array}$$}Equations (a), (b) and (c) follow from \cite[Lemma 2.5]{bz},
Equation (\ref{equa braidedcomulti}) and Lemma \ref{s^2}
respectively. Thus the proof is completed.\qed

From the fact that $\int_{R\#H}^l$ is a right $R\#H$-module, the
following equation holds
$$\xi_R(r)\hdet(h)=\xi_R(h_1(r))\hdet(h_2).$$
We show how $\int^r_{R\#H}$ looks like. Let $\bf{e}'$ be a non-zero
element in $\Ext_R^d(\kk,R)$. There is an algebra homomorphism
$\eta': H\ra\kk$ satisfying $  {\bf e}'\cdot h= \eta'(h)\bf{e}'$ for
all $h\in H$. Applying a similar argument as in the proof of
Proposition \ref{thm int}, we have that if $\int^r_R={}_{\xi'_R}\kk$, then $\int^r_A={\;}_{\xi'}\kk$, where
$\xi'$ is defined by
$\xi'(r\#h)=\xi'_R(\mc{S}^{-1}_H(h_1)(r))\eta'(\mc{S}_H(h_2))$
for all $r\#h\in R\#H$.

Now we give the main theorem of this section.

\begin{thm}\label{cyrtoh}
Let
$H$ be a finite dimensional semisimple  Hopf algebra and  $R$ a Noetherian
braided Hopf algebra in the category $^H_H\mc{YD}$. Suppose that the
algebra $R$ is CY of dimension $d_R$. Then $R\#H$
is CY if and only if the homological determinant  of $R$ is trivial
and the algebra automorphism $\phi$ defined by
$$\phi(r\#h)=\mc{S}_H(r_{(-1)})(\mc{S}_R^2(r_{(0)}))\mc{S}_H^2(h)$$
for all $r\#h\in R\#H$ is an inner automorphism.
\end{thm}

\proof From Proposition \ref{prop cy-as}, we have that $R$ is
AS-regular with $\int_R^l\cong \kk$. In addition, since $H$ is
finite dimensional and semisimple, the algebra $H$ is unimodular.
Thus $\int_H^l=\kk$. Set $A=R\#H$. By Proposition \ref{thm int}, we
obtain that $A$ is AS-regular with $\int_A^l\cong \kk_\xi$, where
$\xi$ is the algebra homomorphism defined by
$\xi(r\#h)=\vps(r)\hdet(h)$ for all $r\#h\in R\#H$. Following
from \cite[Thm. \;2.3]{hvz}, the algebra $A$ is CY if and only if
$\xi=\vps$ and $\mc{S}_A^2$ is an inner automorphism. On one hand,
$\xi=\vps_H$
  if
and only if $\hdet=\vps_H$. On the other hand, by Lemma \ref{s^2},
we have
$\mc{S}_A^2(r\#h)=\mc{S}_H(r_{(-1)})(\mc{S}_R^2(r_{(0)}))\mc{S}_H^2(h)$,
for any $r\#h\in R\#H$.  \qed

One may compare the above theorem with Theorem 2.12 in \cite{liwz}.
Keep the notations as in Theorem \ref{cyrtoh}. When $R$ is
$p$-Koszul and $H$ is involutory, then the automorphism $\phi$ is an
inner automorphism if the homological determinant is trivial. In the
following Example \ref{eg 2}, we  see that the automorphism $\phi$
can be expressed via the homological determinant of the $H$-action.

\begin{eg}\label{eg 2}

 Let $$\mc{D}(\bgm,
(g_i)_{1\se i\se \tt},(\chi_i)_{1\se i\se \tt}, (a_{ij})_{1\se
i,j\se \tt})$$ be a  datum of finite Cartan type (see \cite{as2} for terminology), where $\bgm$ is a
finite abelian group and $(a_{ij})$ is  of type $A_1\times \cdots
\times A_1$. Assume that $V$ is a braided vector space with a basis
$\{x_1,\cdots, x_\tt\}$ whose  braiding is given by
$$c(x_i\ot x_j)=q_{ij}x_j\ot x_i,\;\;1\se i,j\se \tt,$$
where $q_{ij}=\chi_j(g_i)$.

Let $R$ be the following algebra:
$$\kk\lan x_1,\cdots,x_\tt\mid x_ix_j=q_{ij}x_jx_i,\;\;1\se
i<j\se \tt\ran.$$ The algebra $R$ is a braided
Hopf algebra in the category $^\bgm_\bgm\mc{YD}$. Moreover, it is easy to see that $R$ is a Koszul algebra. Assume that
$\mc{K}$ is the Koszul complex (cf. complex (6) in \cite{vdb1})
$$0\ra R\ot R^{!*}_\tt\ra \cdots R\ot R^{!*}_{j}\xra{d_j} R\ot
R^{!*}_{j-1}\cdots\ra R\ot R^{!*}_{1}\ra R\ra 0.$$ Then we have that
$\mc{K}\ra {}_R\kk\ra 0$ is a projective resolution of $\kk$. Each
$R^{!*}_{j}$ is a left $\kk\bgm$-module  with module structure
defined by $$\begin{array}{ccl}[g(\bt)](x^*_{i_1}\wedge\cdots\wedge
x^*_{i_j})&=&\bt(g^{-1}(x^*_ {i_1} \wedge\cdots\wedge x^*_
{i_j}))\\&=&\bt(g^{-1}(x^*_  {i_1} )\wedge\cdots\wedge g^{-1}(x^*_
{i_j} ))\\
&=&(\prod_{t=1}^j\chi_{i_t}(g))\bt(x^*_  {i_1} \wedge\cdots\wedge
x^*_ {i_j} ),\end{array}$$ where $\bt\in S^{!*}_{j} $. Thus each
$R\ot R^{!*}_{j}$ is a left $\kk\bgm$-module. It is not difficult to
see that the differentials in the Koszul complex are also left
$\bgm$-module homomorphisms. By \cite[Prop.
5.0.7]{c}, we have that $\int_R^l\cong R_\tt^{!*}$. Therefore,
$\hdet(g)=\prod_{i=1}^\tt\chi_i(g^{-1})$ for all $g\in \bgm$.

Following from \cite[Prop. 8.2 and Thm. 9.2]{vdb}, the algebra $R$ has the rigid dualizing complex $R_\vph[\tt]$, where $\vph$ is the algebra automorphism defined by $\vph(x_i)=q_{1i}\cdots
q_{(i-1)i}q_{i(i+1)}^{-1}\cdots q_{i\tt}^{-1}x_i$, for $1\se i\se \tt$. By Corollary \ref{cor cyrid}, the algebra $R$ is a CY
algebra if and only if for each $1\se i\se \tt$, $q_{1i}\cdots
q_{(i-1)i}=q_{i(i+1)}\cdots q_{i\tt}$.  In this case, $$\begin{array}{ccl}\hdet(g_j)&=&\prod_{i=1}^\tt\chi_i(g_j^{-1})\\&=&(\prod_{i=1}^{j-1}\chi_i(g_j^{-1}))\chi_j(g_j^{-1})(\prod_{k=j+1}^\tt\chi_k(g_j^{-1}))\\
&=&(\prod_{i=1}^{j-1}q_{ij}) \chi_j(g_j^{-1})(\prod_{k=j+1}^\tt q^{-1}_{jk})\\
&=&\chi_j(g_j^{-1}).\end{array}$$ The algebra automorphism $\phi$ given
in Theorem \ref{cyrtoh} is defined by
$$\phi(x_j)=\chi_j(g_j^{-1})x_j=\hdet(g_j)x_j$$ for all $1\se j\se
\tt$ and $\phi(g)=g$ for all $g\in \bgm$. However, $\chi_j(g_j)\neq
1$ for all $1\se j\se \tt$. The algebra $R\#\kk\bgm$ is not a CY
algebra.
\end{eg}

\begin{eg}\label{cocomm}

Let $\gg$ be a finite dimensional Lie algebra, and $U(\gg)$ the
universal enveloping algebra of $\gg$.  Assume that there is a group
homomorphism $\nu:\bgm\ra Aut_{Lie}(\gg)$, where $Aut_{Lie}(\gg)$ is
the group of Lie algebra automorphisms of $\gg$.  Then it is known
that $U(\gg)\#\kk \bgm$ is a cocommutative Hopf algebra.

It is proved in  \cite[Cor. 3.6]{hvz} that the smash product
$U(\gg)\#\kk \bgm$ is CY if and only if $U(\gg)$ is CY and
$\Im(\nu)\subseteq SL(\gg)$.

Let $d$ be the dimension of $\mk{g}$. By \cite[Lemma 3.1]{hvz}, we have $\int^l_{U(\gg)}\cong \wedge^d
\gg^*$ as left $\bgm$-modules,   where the left $\bgm$-action on
$\mathfrak{g}^*$ is defined by $(g\cdot \al)(x)=\al(g^{-1}x)$ for
all $g\in \bgm$, $\al\in \mathfrak{g}^*$ and $x\in \mathfrak{g}$, and
$\bgm$ acts on $\wedge^d \mathfrak{g}^*$ diagonally. Let
$\{x_1,\cdots,x_d\}$ be a basis of $\gg$.  Then
$$g(x_1^*\wedge\cdots\wedge x_d^*)=\t{det}(\nu(g^{-1}))(x_1^*\wedge\cdots\wedge x_d^*),$$ for all $g\in\bgm$. So
$\hdet(g)=\t{det}(\nu(g))$. That is, if $\Im(\nu)\subseteq SL(\gg)$,
then the homological determinant is trivial.  The algebra
$U(\mk{g})$ is a braided Hopf algebra in the category
$^\bgm_\bgm\mc{YD}$ with trivial coaction. So the automorphism
$\phi$ defined in Theorem \ref{cyrtoh} is the identity. Therefore,
if $U(\mk{g})$ is a CY algebra and $\Im(\nu)\subseteq SL(\gg)$, by
Theorem \ref{cyrtoh}, we can also get that $U(\gg)\#\kk \bgm$ is a
CY algebra.
\end{eg}

\section{Rigid dualizing complexes of braided Hopf algebras over finite group algebras}\label{fgp. fgp}

Let $\bgm$ be a finite group and $R$ a braided Hopf algebra in the category $^\bgm_\bgm\mathcal {YD}$ of Yetter-Drinfeld modules over $\kk\bgm$ such that $R\#\bgm$ is a CY algebra. In this section, we answer the question when the algebra $R$ is a CY algebra.

Let $A$ be a Hopf algebra. By \cite[Appendix, Lemma 11]{pw}, $A$ can
be viewed as a subalgebra of $A^e$ via the algebra homomorphism
$\rho:A\ra A^e$ defined by \begin{equation}\label{rho}\rho(a)=\sum
a_1\ot \mc{S}(a_2).\end{equation} Then $A^e$ is a right   $A$-module
via this embedding. We denote this right $A$-module by
$\mc{R}(A^e)$. Actually, $\mc{R}(A^e)$ is an $A^e$-$A$-bimodule.
Similarly, $A^e$ is also an $A$-$A^e$-bimodule, where the left
$A$-module is induced from the homomorphism $\rho$. Denote this bimodule
by $\mc{L}(A^e)$.

From now on, let $\bgm$ be a finite group and $R$ a braided Hopf
algebra in the category $^\bgm_\bgm\mc{YD}$ with $\bgm$-coaction
$\dt$. The biproduct  $A=R\#\kk\bgm$  is a usual Hopf algebra
\cite{ra}. Let $\mathscr {D}$ be the subalgebra of $A^e$ generated
by the elements of the form $(r\#g)\ot (s\#g^{-1})$ with $r,s\in R$
and $g\in \bgm$.


Note that  $R$ is a $\bgm$-graded module, i.e., $R =\op_{g\in
\bgm}R_g$, where $R_g = \{r \in R \mid \dt(r) = g \ot r\}$.
Therefore, for any $r\in R$, it can be written as  $r=\sum_{g\in
\bgm}r_g$ with $r_g\in R_g$. Then $\dt(r)=\sum_{g\in \bgm}g\ot r_g$.

\begin{lem}
The subalgebra $\mathscr {D}$ is a left and right $A$-submodule of
$\mc{L}(A^e)$ and $\mc{R}(A^e)$ respectively.
\end{lem}
\proof For any $r\#h\in A$, by equations (\ref{equa braidedcomulti})
and (\ref{equa braidedanti}), we have $$\Delta(r\#h)=\sum_{g\in
\bgm} r^1\#gh\ot (r^2)_g\#h$$ and $$\mc{S}_A(r\#h)=\sum_{g\in \bgm}
h^{-1}g^{-1} \mc{S}_R(r_g).$$ Any element in $\ms{D}$ can be written
as a linear combination of elements of the form  $s\#k\ot
t\#k^{-1}\in \ms {D}$ with $s,t\in R$ and $k\in \bgm$.
$$\begin{array}{cl}&(r\#h)\cdot(s\#k\ot
t\#k^{-1})\\=&\sum_{g\in \bgm}(r^1\#gh)(s\#k)\ot (t\#k^{-1})\mc{S}_A((r^2)_g\#h)\\
=&\sum_{g\in \bgm}(r^1\#gh)(s\#k)\ot (t\#k^{-1})h^{-1}g^{-1}\mc{S}_R((r^2)_g)\\
=&\sum_{g\in \bgm}(r^1{(gh)(s)}\#ghk)\ot (t(k^{-1}h^{-1}g^{-1})(\mc{S}_R((r^2)_g))\#k^{-1}h^{-1}g^{-1})\\
\in &\mathscr {D}.  \end{array}$$ This shows that $\mathscr{D}$ is a
left $A$-submodule of $\mc{L}(A^e)$. Similarly,  $\mathscr{D}$ is
also a right $A$-submodule of $\mc{R}(A^e)$. \qed

The following lemma is known, we include it for completeness.
\begin{lem}\label{lem d}
\begin{enumerate}
\item[(a)] Both $\mc{L}(A^e)$ and $\mc{R}(A^e)$ are free as $A$-modules.
\item[(b)] $\mc{R}(A^e)\ot_A \kk\cong A$ as left $A^e$-modules and this isomorphism restricts to the isomorphism $\ms{D}\ot_A\kk\cong R$.
\item[(c)] If $\xi:A\ra\kk$ is  an algebra homomorphism, then there is an isomorphism $\kk_\xi\ot_A \mc{L}(A^e)\cong {A_{[\xi]\mc{S}_A^2}}$ of right $A^e$-modules and the isomorphism restricts to an isomorphism
$\kk_\xi\ot_A \ms{D}\cong {R_{([\xi]\mc{S}_A^2)|_R}}$.
\end{enumerate}
\end{lem}
\proof (a) was proved in \cite[Lemma 2.2]{bz}. The module $L(A^e)$
defined in the same paper is isomorphic to $\mc{R}(A^e)$ as right
$A$-modules. It was proved  that $\varphi:A_A\ot A^{op}\ra
\mc{R}(A^e)$ defined by $\vph(a\ot b)=a_1\ot b\star\mc{S}_A(a_2)$ is
an isomorphism, where $\star$ denotes the multiplication in
$A^{op}$. The right $A$-module structure on $A_A\ot A^{op}$ is
defined by $(a\ot b)\cdot c=ac\ot b$ for all $a$, $b$ and $c\in A$.
Similarly, $\mc{L}(A^e)\cong {_AA}\ot A^{op}$ as free left
$A$-module.

(b) The isomorphism  $\mc{R}(A^e)\ot_A \kk\cong A$ of left $A^e$-modules can be
found in \cite[Appendix, Lemma 11]{pw}.  The homomorphism
$\phi:\mc{R}(A^e)\ot_A \kk\ra A$ given by $\phi(a\ot b \ot 1)=ab$ is
an $A^e$-isomorphism. It is clear that $\psi$ restricts to an
isomorphism from $\ms{D}\ot_A\kk$ to $R$.

(c) It was proved in   \cite[Lemma 4.5]{bz} that $\kk_\xi\ot_A
\mc{L}(A^e)\cong {A_{[\xi]\mc{S}_A^2}}$ as right $A^e$-modules.
Here we give another proof. We construct the  the isomorphism
explicitly. Define a homomorphism $\Phi:\kk_\xi\ot_A A^e\ra
{A_{[\xi]\mc{S}_A^2}}$ by $\Phi(1\ot a\ot
b)=\xi(a_1)b\mc{S}_A^2(a_2)$ and a homomorphism
$\Psi:{A_{[\xi]\mc{S}_A^2}}\ra \kk_\xi\ot_A A^e$ by $\Psi(a)=1\ot 1\ot a$. For any $x,a, b\in A$, we have $$\begin{array}{ccl} \Phi(1\ot x_1a\ot b\mc{S}(x_2))&=& \xi(x_1)\xi(a_1)\ot b\mc{S}(x_3)\mc{S}^2(x_2)\mc{S}^2(a_2)\\
&=& \xi(x_1)\xi(a_1)\ot b\mc{S}(\varepsilon(x_2))\mc{S}^2(a_2)\\
&=& \xi(x)\xi(a_1)\ot b\mc{S}^2(a_2)\\
&=&\xi(x)\Phi(1\ot a\ot b).\end{array}$$
This shows that $\Phi$ is
well defined. Similar calculations show that $\Phi$ and $\Psi$ are
right $A^e$-module homomorphisms and they are inverse to each other.

It is straightforward to check that the isomorphism $\kk_\xi\ot_A
\mc{L}(A^e)\cong {A_{[\xi]\mc{S}_A^2}}$ restricts to the isomorphism
$\kk_\xi\ot_A \ms{D}\cong {R_{([\xi]\mc{S}_A^2)|_R}}$. \qed

\begin{lem}\label{lem fgp}
$\Hom_{R^e}(\mathscr{D}, R^e)\cong \mathscr{D}$ as
$A$-$R^e$-bimodules.
\end{lem}
\proof The algebra $\mathscr{D}$ is   an $A$-$R^e$-bimodule. Note
that the $A$-module structure is  induced from the homomorphism
$\rho$ defined in (\ref{rho}). On the other hand, the $A$-$R^e$-bimodule
structure on $\Hom_{R^e}(\mathscr{D}, R^e)$ is induced from the
right $A$-module structure on $\mathscr{D}$ and the right
$R^e$-module structure on $R^e$.  We have
$r\#g=(1\#g)(g^{-1}(r)\#1)$ for any $r\#g \in R\#\kk\bgm$.
Therefore, an element in $\mathscr{D}$ can be expressed of the form
$\sum_{g\in \bgm} (1\#g^{-1})(r^g\#1)\ot s^g \#g$ with $r^g, s^g\in
R$. For simplicity, we write $gr$ for  the element $(1\#g)(r\#1)$ with $r\in
R$ and $g\in \bgm$. Let $\Psi:\mathscr{D}\ra
\Hom_{R^e}(\mathscr{D}, R^e)$ be a homomorphism defined by
$$[\Psi(\sum_{g\in \bgm} g^{-1}r^g\ot (s^g \#g))](h\ot h^{-1})=r^h\ot s^h,$$ for
$\sum_{g\in\bgm} g^{-1}r^g\ot s^g \#g\in \mathscr{D}, h\ot h^{-1}\in
\mathscr{D}$. Next define a homomorphism $\Phi:
\Hom_{R^e}(\mathscr{D}, R^e)\ra\mathscr{D}$  as follows:
 $$\Phi(f)=\sum_{g\in
\bgm}( g^{-1}\ot g)f(g\ot g^{-1})$$ for $f\in
\Hom_{R^e}(\mathscr{D}, R^e)$. It is clear that $\Phi$ is a right
$R^e$-homomorphism. On the other hand, we have
 {\small $$\begin{array}{ccl}\Phi((r\#h)f)&=&\sum_{g\in \bgm}(g^{-1}\ot g)((r\#h)f)(g\ot g^{-1})\\
&=&\sum_{g\in \bgm}\sum_{k\in \bgm}(g^{-1}\ot g)f(g(r^1\#k)h\ot \mc{S}_A((r^2)_k\#h)g^{-1})\\
&=&\sum_{g\in \bgm}\sum_{k\in \bgm}(g^{-1}\ot g)f(g(r^1\#k)h\ot h^{-1}k^{-1}\mc{S}_R((r^2)_k)g^{-1})\\
&=&\sum_{g\in \bgm}\sum_{k\in \bgm}(g^{-1}\ot g)f(g(r^1)\#gkh\ot
h^{-1}k^{-1}g^{-1}{g(\mc{S}_R((r^2)_k))}),
\end{array}$$}
and
{\small $$\begin{array}{cl}&(r\#h)\Phi(f)\\=&(\sum_{k\in \bgm}r_1\#kh\ot h^{-1}k^{-1}\mc{S}_{R}(r_{2k}))\sum_{g\in \bgm} (g^{-1}\ot g)f(g\ot g^{-1})\\
=&\sum_{k\in \bgm}\sum_{g\in \bgm}(r^1\#khg^{-1}\ot gh^{-1}k^{-1}\mc{S}_{R}((r^2)_k))f(g\ot g^{-1})\\
=&\sum_{k\in \bgm}\sum_{g\in \bgm}(khg^{-1}{(gh^{-1}k^{-1})(r^1)}\ot (gh^{-1}k^{-1})\mc{S}_{R}((r^2)_k)gh^{-1}k^{-1})f(g\ot g^{-1})\\
=&\sum_{k\in \bgm}\sum_{g\in \bgm}(khg^{-1}\ot
gh^{-1}k^{-1})f((gh^{-1}k^{-1})(r^1)\#g\ot
g^{-1}(gh^{-1}k^{-1})(\mc{S}_{R}((r^2)_k)))\\
=&\sum_{g\in \bgm}\sum_{k\in \bgm}(g^{-1}\ot g)f(g(r^1)\#gkh\ot
h^{-1}k^{-1}g^{-1}{g(\mc{S}_R((r^2)_k))}).
\end{array}$$}
So $\Phi$ is an $A$-$R^e$-bimodule homomorphism. It is easy to check that
$\Phi$ and $\Psi$ are inverse to each other. Thus $\Phi$ is an
isomorphism. \qed

\begin{lem}\label{bimod}
Let $\bgm$ be a finite group and $R$  a braided Hopf algebra in the
category $^\bgm_\bgm\mc{YD}$. If $A=R\#\kk\bgm$ is AS-Gorenstein
with $\int_A^l\cong\kk_\xi$, where $\xi:A\ra\kk$ is an algebra homomorphism, then we have $R$-$R$-bimodule
isomorphisms
$$ \Ext^i_{R^e}(R,R^e)\cong\begin{cases}0,&i\neq d;\\
{R_{([\xi]\mc{S}_A^2)|_R}},&i=d.\end{cases}$$
\end{lem}
\proof  We have the following isomorphisms,
$$\begin{array}{ccl} \Ext^i_{R^e}(R,R^e)&\cong&\Ext^i_{R^e}(\ms{D}\ot_A
\kk,R^e)\\&\cong&\Ext^i_A({}_A\kk,\Hom_{R^e}(\ms{D} ,R^e))\\
&\cong&\Ext^i_A({}_A\kk,\ms{D})\\
&\cong&\Ext^i_A({}_A\kk,A)\ot_A\ms{D}\\
&\cong&\begin{cases}0,&i\neq d;\\\kk_\xi\ot_A \ms{D}\cong
{R_{([\xi]\mc{S}_A^2)|_R}},&i=d.\end{cases}\end{array}$$ The first, the third
and the last isomorphism follow from Lemma \ref{lem d},  Lemma \ref{lem
fgp} and  Lemma \ref{lem d} respectively.\qed

\begin{lem}\label{hs}
If the projective dimension of $A=R\#\kk\bgm$ is finite and $R$ is
Noetherian, then $R$ is homologically smooth.
\end{lem}
\proof By assumption, the algebra $A$ is Noetherian, and $_A\kk$ has
a finite projective resolution
$$0\ra P_d\ra P_{d-1}\ra\cdots\ra P_1\ra P_0\ra \kk\ra 0,$$ such that each $P_i$, $0\se i\se d$, is a finitely generated projective $A$-module.  $\ms{D}$ is a summand of $A^e$ as a right $A$-module. Indeed, $A^e\cong \op_{h\in \bgm}\ms{D}^h$, where $\ms{D}^h$ is the right  $A$-submodule of $A^e$ generated by elements of the form $(r\#hg)\ot (s\#g^{-1})$ with  $r,s\in R$ and $g\in \bgm$.  By Lemma \ref{lem d}, $A^e$ is free as a right $A$-module. Therefore, the functor $\ms{D}\ot_A- $ is exact. We obtain an exact sequence \begin{equation}\label{comp3}0\ra \ms{D}\ot_AP_d\ra \ms{D}\ot_AP_{d-1}\ra\cdots\ra \ms{D}\ot_AP_1\ra \ms{D}\ot_AP_0\ra \ms{D}\ot_A\kk\ra 0.\end{equation} $\ms{D}$ is projective as left $R^e$-module and $\ms{D}\ot_A\kk\cong R$ as  left $R^e$-modules (Lemma \ref{lem d}). So the complex (\ref{comp3}) is a projective bimodule resolution of $R$. Because each $P_i$ is a finitely generated $A$-module and $\bgm$ is a finite group, each $\ms{D}\ot_AP_{i}$ is a finitely generated   left $R^e$-module. Therefore, we conclude that $R$ is homologically smooth. \qed

The homological integral  of the  skew group algebra $R\#\kk\bgm$
was discussed by He, Van Oystaeyen and Zhang in \cite{hvz}. Based on
their work,  we use the homological determinant of the group
action to describe the homological integral of $R\#\kk\bgm$.

\begin{lem}\label{int fgp}
Let $\bgm$ be a finite group and $R$  a braided Hopf algebra in the
category $^\bgm_\bgm\mc{YD}$. If  $R$ is an AS-Gorenstein algebra
with injective dimension $d$ and $\int_R^l\cong \kk_{\xi_R}$, where
$\xi_R:R\ra \kk$ is an algebra homomorphism,  then the algebra
$A=R\#\kk\bgm$ is AS-Gorenstein with injective dimension $d$ as
well, and $\int_A^l\cong \kk_\xi$, where $\xi:A\ra \kk$ is the algebra homomorphism  defined by
$\xi(r\#h)=\xi_R(r)\hdet(h)$ for any $r\#h\in R\#\kk\bgm$.
\end{lem}
\proof By  \cite[Prop. 1.1 and 1.3]{hvz},  we have that
$A=R\#\kk\bgm$ is AS-Gorenstein of injective dimension $d$,
$\int_R^l$ is a
1-dimensional left $\bgm$-module, and as right $A$-modules:\vspace{3mm}\\
\centerline{$\int^l_A\cong(\int_R^l\ot\kk\bgm)^\bgm$,}\vspace{3mm}\\ where the
right $A$-module structure on $\int_R^l\ot\kk\bgm$ is defined by $$({  e}\ot
g)\cdot(r\#h)={  e}(g(r))\ot gh,$$ for $g\in \kk\bgm$,
$r\#h\in R\#\kk\bgm$ and ${e}\in \int_R^l$, and the left $\bgm$-action on $\int_R^l\ot\kk\bgm$ is the diagonal
one. Let ${\bf e}$ be a basis of $\int_R^l$. It is not difficult to check
that the element  $\sum_{g\in \bgm}{g({\bf
e})}\ot g$  is a basis of $(\int_R^l\#\kk\bgm)^\bgm$. Let $\eta:\kk
\bgm\ra\kk$ be an algebra homomorphism such that $h\cdot {\bf
e}=\eta(h){\bf e}$ for all $h\in\bgm$.
For any $r\#h\in R\#\kk\bgm$, we have $$\begin{array}{ccl}(\sum_{g\in \bgm}{g({\bf e})}\#g)(r\#h)&=&\sum_{g\in \bgm}{g({\bf e})}g(r)\#gh\\
&=&\sum_{g\in \bgm}{g({\bf e}r)}\#gh\\
&=&\xi_R(r)\sum_{g\in \bgm}{g({\bf e})}\#gh\\
&=&\xi_R(r)\eta(h^{-1})\sum_{g\in \bgm}{(gh)({\bf e})}\#gh\\
&=&\xi_R(r)\eta(h^{-1})\sum_{g\in \bgm}{g({\bf e})}\#g\\
&=&\xi_R(r)\hdet(h)\sum_{g\in \bgm}{g({\bf e})}\#g\\
&=&\xi(r\#h)\sum_{g\in \bgm}{g({\bf e})}\#g.
\end{array}$$
It implies that $\int_A^l\cong\kk_\xi$. \qed

The following proposition shows that the AS-regularity of $R\#\kk\bgm$ depends strongly
on the AS-regularity of $R$.

\begin{prop}\label{asreg}
Let $\bgm$ be a finite group and $R$ a braided Hopf algebra in the
category $_\bgm^\bgm\mc{YD}$. Then  $R$ is AS-regular if and only if
$A=R\#\kk\bgm$ is AS-regular.
\end{prop}
\proof Assume that $R$ is AS-regular. By Lemma \ref{int fgp}, the algebra
$A$ is AS-Gorenstein. To show that $A$ is AS-regular, it suffices to
show that the global dimension of $A$ is finite.  Since the global
dimension of $R$ is finite, there is a finite projective resolution
of $\kk$ over $R$,
$$0\ra P_d\ra P_{d-1}\ra\cdots P_1\ra P_0\ra \kk\ra 0.$$ Note that $A$ is projective as a right $R$-module. Tensoring this resolution with $A\ot_R-$, we obtain an exact sequence
$$0\ra A\ot_R P_d\ra A\ot_R P_{d-1}\ra\cdots A\ot_R P_1\ra A\ot_R P_0\ra A\ot_R \kk\ra 0.$$ It is clear that each $A\ot_R P_i$ is projective. This shows that the projective dimension of $A\ot_R\kk$ is finite. But $_A\kk$ is a direct summand of $A\ot_R\kk$ as an $A$-module (\cite[Lemma III.4.8]{ars}). So the projective dimension of $_A\kk$ is finite. Since $A$ is a Hopf algebra, the global dimension of $A$ is finite.

Conversely, if $A$ is AS-regular, then $R$ is AS-regular by  Lemma \ref{bimod}, Lemma
\ref{hs} and Remark \ref{twisted}. \qed

We are ready to give the rigid dualizing complex of an AS-Gorenstein braided Hopf
algebra.

\begin{thm}\label{thm ridfgp}
Let $\bgm$ be a finite group and $R$  a braided Hopf algebra in the
category $^\bgm_\bgm\mc{YD}$. Assume that $R$ is an
AS-Gorenstein algebra with  injective dimension $d$. If
$\int_R^l\cong\kk_{\xi_R}$, for some algebra homomorphism
$\xi_R:R\ra\kk$,  then $R$ has a  rigid dualizing complex $_\vph
R[d]$, where $\vph$ is the algebra automorphism defined by $$\vph(r)=\sum_{g\in
\bgm}\xi_R(r^1)\hdet(g){{g^{-1}}(\mc{S}_R^2((r^2)_g))},$$ for any $r\in R$.
\end{thm}
\proof Let $A$ be $R\#\kk\bgm$.  By Lemma \ref{int fgp}, $A$ is AS-Gorenstein with $\int_A^l\cong\kk_\xi,$ where $\xi:A\ra \kk$ is the algebra homomorphism
defined by $$\xi(r\#h)=\xi_R(r)\hdet(h)$$ for any $r\#h\in R\#\kk\bgm$. By Lemma \ref{bimod}, there are bimodule isomorphisms $$ \Ext^i_{R^e}(R,R^e)\cong\begin{cases}0,&i\neq d;\\
R_{([\xi]\mc{S}_A^2)|_R},&i=d.\end{cases}$$
For any $r\in R$,
$$\begin{array}{ccl}[\xi]\mc{S}_A^2(r)&=&\sum_{g\in \bgm}\xi(r^1\#g)\mc{S}_A^2((r^2)_g)\\&=&\sum_{g\in \bgm}\xi_R(r^1)\hdet(g)\mc{S}_A^2((r^2)_g)\\
&=&\sum_{g\in
\bgm}\xi_R(r^1)\hdet(g)g^{-1}(\mc{S}_R^2((r^2)_g)).\end{array}$$ Now
the theorem follows from Lemma \ref{vdb}. \qed

\begin{rem}
The algebra $A=R\#\kk\bgm$ has a rigid dualizing complex
$_{[\xi]\mc{S}^2_A}A[d]$ (\cite[Prop. \;4.5]{bz}). Observe that
the algebra automorphism $\vph$  given in Theorem \ref{thm ridfgp} is just
the restriction of $[\xi]\mc{S}_A^2$ on $R$.
\end{rem}

Now we can characterize the CY property of $R$ in case $R\#\kk\bgm$
is CY.

\begin{thm}\label{thm cyator}
Let $\bgm$ be a finite group and $R$  a Noetherian braided Hopf algebra in the
category $^\bgm_\bgm\mc{YD}$. Define an algebra automorphism $\vph$
of $R$ by $$\vph(r)=\sum_{g\in \bgm}g^{-1}(\mc{S}_R^2(r_g)),$$ for any
$r\in R$. If $R\#\kk\bgm$ is a CY algebra, then $R$ is CY if and
only if the algebra automorphism $\vph$ is an inner automorphism.
\end{thm}

\proof Assume that $A=R\#\kk\bgm$ is a CY algebra of dimension $d$. By
Proposition \ref{prop cy-as}, $A$ is AS-regular of global dimension
$d$ and $\int_A^l\cong \kk$. It follows from Lemma \ref{hs} that $R$
is homologically smooth.

Since $\int^l_A\cong \kk$, by Lemma \ref{bimod}   there are
$R$-$R$-bimodule isomorphisms
$$ \Ext^i_{R^e}(R,R^e)\cong\begin{cases}0,&i\neq d;\\
{R_{ \mc{S}_A^2|_R}},&i=d.\end{cases}$$ Following Remark
\ref{twisted}, we obtain that $R$ is AS-regular. Suppose
$\int^l_R\cong \kk_{\xi_R}$ for some algebra homomorphism
$\xi_R:R\ra \kk$. Then by Lemma \ref{int fgp}, $\int^l_A\cong
\kk_\xi$, where $\xi:A\ra\kk$ is defined by
$\xi(r\#h)=\xi_R(r)\hdet(h)$ for any $r\#h\in R\#\kk\bgm$. But
$\int^l_A\cong \kk$. Therefore, $\xi_R=\vps_R$ and
$\t{hdet}=\vps_H$. It follows from Theorem \ref{thm ridfgp} that the
rigid dualizing complex of $R$ is isomorphic to $_\vph R[d]$, where
$\vph$ is defined by
$$\begin{array}{ccl}\vph(r)&=&\sum_{g\in
\bgm}\xi_R(r^1)\hdet(g){{g^{-1}}(\mc{S}_R^2((r^2)_g))}\\
&=&\sum_{g\in \bgm}g^{-1}(\mc{S}_R^2(r_g))\end{array}$$ for any $r\in
R$. Now the theorem follows from Corollary \ref{cor cyrid}. \qed

\begin{cor}\label{cor iff}
Let $\bgm$ be a finite group and $R$  a braided Hopf algebra in the
category $^\bgm_\bgm\mc{YD}$. Assume that $R$ is an
AS-regular algebra. Then the following two conditions are
equivalent:
\begin{enumerate}
\item[(a)] Both $R$ and $R\#\kk\bgm$ are CY algebras.
\item[(b)] The following three conditions are satisfied:
\begin{enumerate}
\item[(i)] $\int_R^l\cong \kk$;
\item[(ii)] The homological determinant of the group action is
trivial;
\item[(iii)]  The algebra automorphism $\vph$ defined by $$\vph(r)=\sum_{g\in \bgm}g^{-1}(\mc{S}_R^2(r_g))$$ for all $r\in R$ is
an inner automorphism.
\end{enumerate}
\end{enumerate}
\end{cor}
\proof $(a)\Rightarrow(b)$ Since $R$ is a CY algebra, by Proposition \ref{prop cy-as} we have $\int_R^l\cong \kk$.  Because both $R$ and $R\#\kk\bgm$ are CY, (ii) and (iii) are satisfied by Theorem \ref{cyrtoh} and Theorem \ref{thm cyator}.

$(b)\Rightarrow(a)$ Since $R$ is AS-regular, $R\#\kk\bgm$ is
AS-regular by Proposition \ref{asreg}. Thus $R$ is homologically
smooth (Lemma \ref{hs}). By Theorem \ref{thm ridfgp}, if the three
conditions in (b) are satisfied, then the rigid dualizing complex of
$R$ is isomorphic to $R[d]$, where $d$ is the injective dimension of
$R$. So $R$ is a CY algebra. That the algebra $R\#\kk\bgm$ is a CY
algebra follows from Theorem \ref{cyrtoh}.  \qed

\begin{eg}
Keep the notations from Example \ref{cocomm}. Assume that $\bgm$
is a finite group and that $\mk{g}$ is a finite dimensional $\bgm$-module
Lie algebra. Suppose there is a group homomorphism $\nu:\bgm\ra
Aut_{Lie}(\mk{g})$.  In Example \ref{cocomm}, we use Theorem
\ref{cyrtoh} to obtain that if $U(\mk{g})$ is a CY algebra and
$\Im(\nu)\subseteq SL(\mk{g})$ then $U(\mk{g})\#\kk\bgm$ is a CY
algebra. Now by Theorem \ref{thm cyator}, if $U(\mk{g})\#\kk\bgm$ is
a CY algebra, then $U(\mk{g})$ is a CY algebra as well. This is
because $U(\mk{g})$ is a braided Hopf algebra in $^\bgm_\bgm\mc{YD}$
with trivial coaction,  the algebra automorphism $\vph$ in Theorem
\ref{thm cyator} is the identity.

By \cite[Prop. 6.3]{bz}, we have that $\int_{U(\mk{g})}^l=\kk_\xi$,
where $\xi(x)=\t{tr}(\ad(x))$ for all $x\in \mk{g}$. We calculate in
Example \ref{cocomm} that $\hdet(g)=\t{det}(\nu(g))$ for $g\in
\bgm$.   Therefore, both $U(\mk{g})$ and $U(\mk{g})\#\kk\bgm$ are CY
algebras if and only if $\t{tr}(\ad(x))=0$ for all $x\in \mk{g}$ and
$\Im(\nu)\subseteq SL(\mk{g})$. This coincides with Corollary 3.5 and
Lemma 4.1 in \cite{hvz}.

\end{eg}

We refer to \cite{as2} for the definition of a datum of finite Cartan type
and the definition of the algebras $U(\mc{D},\lmd)$. The algebras $U(\mc{D},\lmd)$
were constructed to classify finite dimensional pointed Hopf algebras whose group-like
elements form  an abelian group.

Let $$\mc{D}(\bgm, (g_i)_{1\se i\se \tt},(\chi_i)_{1\se i\se \tt},
(a_{ij})_{1\se i,j\se \tt})$$ be a datum of finite Cartan type for a
finite abelian group $\bgm$ and $\lmd$ a family of linking parameters for $\mc{D}$. Let $\{\al_1,\cdots,\al_\tt\}$ be a set
of simple roots  of the root system corresponding to the  Cartan
matrix $(a_{ij})$. Assume that $w_0=s_{i_1}\cdots s_{i_p}$ is a reduced
decomposition of the longest element in the Weyl group as a product
of simple reflections. Then the positive roots are as follows
$$\bt_1=\al_{i_1},\bt_2=s_{i_1}(\al_{i_2}),\cdots, \bt_p=s_{i_1}\cdots s_{i_{p-1}}(\al_{i_p}).$$
If
$\bt_i=\sum_{i=1}^{\tt} m_i\al_i$, then we define
$\chi_{_{\bt_i}}={\chi}_1^{m_1}\cdots {\chi}_\tt^{m_\tt}$.

The following proposition characterizes the CY property of the algebra $U(\mc{D},\lmd)$.

\begin{prop}
(a) The algebra $A=U(\mc{D},\lmd)$ is
AS-regular of global dimension $p$ and $\int^l_A=\kk_\xi$, where
$\xi$ is the algebra homomorphism defined by
$\xi(g)=(\prod_{i=1}^p\chi_{_{\bt_i}})(g)$, for all $g\in \bgm$ and
$\xi(x_i)=0$ for all $1\se i\se \tt$.

(b) The algebra $A$ is CY if and only if
$\prod_{i=1}^p\chi_{_{\bt_i}}=\vps$  and $\mc{S}_A^2$ is an inner
automorphism.
\end{prop}

\proof (a) can be obtained by applying \cite[Thm. 2.6]{as2} and  a similar argument as in the
proof of Theorem 2.2 in \cite{yuz}. (b) follows from \cite[Thm. 2.3]{hvz}. \qed


Let $R$ be the algebra generated by $x_1,\cdots ,x_\tt$ subject to
the relations
$$(\ad_cx_i)^{1-a_{ij}}(x_j)=0,1\se i,j\se \tt,\;\;i\neq j.$$ Then $U(\mc{D},0)=R\#\kk\bgm$, where $U(\mc{D},0)$ is the associated graded algebra of $U(\mc{D},\lmd)$ with respect to its coradical filtration.

\begin{prop}
The algebra $R$ is CY if and only if $\prod_{i=1,i\neq
j_k}^p\chi_{_{\bt_i}}(g_k)=1$ for each $1\se k\se \tt$.
\end{prop}

\proof By Lemma \ref{hs} and Theorem \ref{thm ridfgp}, we have that $R$ is homologically smooth, and that it has a rigid dualizing complex $_\vph R[p]$, where $\vph$ is the restriction of $[\xi]\mc{S}_A^2$ on $R$. That is, $\vph$ is defined by $\vph(x_k)=\prod_{i=1,i\neq j_k}^p\chi_{_{\bt_i}}(g_k)(x_k)$, $1\se k\se \tt$,   where each $1\se
j_k\se p$   is the integer such that $\bt_{j_k}=\al_k$. Therefore,
$R$ is CY if and only if $\prod_{i=1,i\neq
j_k}^p\chi_{_{\bt_i}}(g_k)=1$ for each $1\se k\se \tt$.\qed

One may compare these results with Theorem 2.3, Theorem
3.9 and Lemma 4.1 in \cite{yuz}.

In case the algebra $U(\mc{D},0)=R\#\kk\bgm$ is CY, the algebra automorphism $\vph$ defined in Theorem \ref{thm cyator} is $\vph(x_i)=\chi_i(g_i^{-1})(x_i)$, $1\se i\se \tt$. However, $\chi_i(g_i)\neq 1$ for all $1\se i\se \tt$. We conclude that when $R\#\kk\bgm$ is CY, the algebra $R$ is not a CY algebra.

Now we give two examples of CY pointed Hopf algebra with a finite group of group-like elements.
\begin{eg} Let $A$ be $U(\mc{D},\lmd)$ with
the datum $(\mc{D},\lmd)$   given by
\begin{itemize}
\item $\bgm=\lan y_1,y_2 \ran\cong \ZZ_2\times \ZZ_2$;
\item The Cartan matrix is of type $A_2$;
\item $g_i=y_i$, $1\se i\se 2$;
\item $\chi_i$, $1\se i\se 2$, are given by the following table. $$\begin{tabular}{|c|c|c|}\hline
&$y_1$&$y_2$\\\hline $\chi_1$&$-1$&$1$\\\hline
$\chi_2$&$-1$&$-1$\\\hline
\end{tabular}$$
\item $\lmd=0$
\end{itemize}
The algebra $A$ is a CY algebra of dimension 3.
Let $R$ be the algebra generated by $x_1$ and $x_2$
subject to relations
$$x_1^2x_2-x_2x_1^2=0\t{ and }
x_2^2x_1-x_1x_2^2=0.$$ Then $A=R\#\kk\bgm$.
The rigid dualizing complex of $R$ is $_\vph R[3]$, where
$\vph=-\id$.
\begin{rem}From the proof of Proposition 5.8 in \cite{yuz}, we can see that if $A=U(\mc{D},\lmd)$ is a CY algebra such that  $(\mc{D},\lmd)$ is a generic datum, then the Cartan matrix in $\mc{D}$ cannot be of type $A_2$.
\end{rem}

\end{eg}

\begin{eg} Let $A$ be $U(\mc{D},\lmd)$ with
the datum $(\mc{D},\lmd)$   given by
\begin{itemize}
\item $\bgm=\lan y_1,y_2 \ran\cong \ZZ_n\times \ZZ_n$;
\item The Cartan matrix is of type $A_1\times A_1$;
\item $g_i=y_i$, $i=1,2$;
\item $\chi_1(y_i)=q$,  $\chi_2(y_i)=q^{-1}$, $i=1,2$, where $q\in \kk$ is an $n$-th root of unity;
\item $\lmd=1$.
\end{itemize}
The algebra $A$ is a CY algebra of dimension 2.

Let $R$ be the algebra
$\kk\lan x_1,x_2\mid x_1x_2=q^{-1}x_2x_1\ran.$ Then $\Gr A=U(\mc{D},0)=R\#\kk\bgm$, where $\Gr A$ is the associated graded algebra of $A$ with respect to the coradical filtration of $A$. The rigid dualizing complex of $R$ is $_\vph R[3]$, where $\vph$ is defined by $\vph(x_1)=q^{-1}x_1$ and $\vph(x_2)=qx_2$.
\end{eg}

\vspace{5mm}

\subsection*{Acknowledgement} This work forms a part of the  PhD thesis of the first author at Hasselt University. She would like to thank BOF of UHasselt its the financial support.

\vspace{5mm}

\bibliography{}

\end{document}